\documentstyle[12pt]{article}

  \newcommand{\const}{\rm const}

\begin{document}

   \begin{center}

 {\bf Restricted version of Grand Lebesgue Spaces }

\vspace{5mm}

 {\bf M.R.Formica, \ E.Ostrovsky, \ L.Sirota }

\end{center}

\vspace{6mm}

 \ Universit\`{a} degli Studi di Napoli Parthenope, via Generale Parisi 13, Palazzo Pacanowsky, 80132,
Napoli, Italy \\

e-mail: mara.formica@uniparthenope.it \\

\vspace{3mm}

 \ Israel,  Bar-Ilan University, department  of Mathematic and Statistics, 59200, \\

\vspace{3mm}

e-mails: \\
 eugostrovsky@list.ru \\
sirota3@bezeqint.net \\

\vspace{6mm}

 \begin{center}

 \ {\bf Abstract}

\vspace{4mm}

\end{center}

 \ We introduce a so-called restricted, in particular, discrete version of (Banach) Grand Lebesgue Spaces (GLS), investigate
 its properties and derive the conditions of coincidence with the classical ones.\par
  \ We show also that these spaces forms also a Banach algebra relative the convolution operation on the unimodular local compact topological
 group equipped with Haar's measure, alike in the complete GLS case.\par

\vspace{5mm}

\begin{center}

 \ {\it Key words and phrases}\

\vspace{4mm}

\end{center}

 \ Measurable and probability spaces, random variable (r.v.), Lebesgue (Riesz) norm, classical and discrete Grand Lebesgue Spaces (GLS), slowly varying function,
 partition, bilateral estimations, comparison and equivalence of norms, unimodular local compact topological group, convolution, Haar's measure, Banach algebra.\par

\vspace{5mm}

\section{Definitions. Statement of problem.}

\vspace{5mm}

 \ Let $ \ (\Omega = (\omega), F,{\bf P})  \  $ be non-trivial probability space, and let $ \ \psi = \psi(p), \ 1 \le p < \infty \ $ be  finite continuous
{\it strictly increasing} positive numerical valued function  such that

\begin{equation} \label{positivity}
 \psi(1) = \inf_{p \in [1,\infty) } \psi(p)  = 1.
\end{equation}

 \ The set of all such a functions will be denote by $ \ \Psi: \ \Psi = \{\psi(\cdot)\}. \ $ \par

 \ The ordinary Lebesgue-Riesz  norm $ \ |\xi|_p \ $ of a r.v. $ \ \xi \ $ will be defined as usually

$$
|\xi|_p := \left[ \ {\bf E}|\xi|^p \  \right]^{1/p}, \ 1 \le p < \infty.
$$

 \ Recall the classical definition of the  Grand Lebesgue Spaces (GLS).
 By definition, the random variable (function)
 $ \ \xi: \ \Omega \to R \ $ belongs to the space $ \ G\psi,  \ $ iff it has a finite GLS norm

\begin{equation} \label{GLS norm}
||\xi||G\psi \stackrel{def}{=} \sup_{p \ge 1} \left\{ \ \frac{|\xi|_p}{\psi(p)} \ \right\}.
\end{equation}

 \ The function $\psi =\psi(p)$ is named ordinary  as  a {\it generating
function} for the Grand Lebesgue Space $G\psi$. \par

\vspace{4mm}

  \ These GLS spaces are rearrangement-invariant Banach functional spaces in
the classical sense, \cite{Bennett Sharpley1988}, chapters 1,2;  \ and were investigated in particular in  many
works, see e.g. \cite{Buld Mush OsPuch1992}, \cite{caponeformicagiovanonlanal2013}, \cite{Fiorenza1}, \cite{Fiorenza2},
\cite{fioforgogakoparakoNAtoappear}, \cite{Fiorenza2018a}, \cite{Fiorenza2018b}, \cite{Formica As or not}, \cite{Iwaniec1},
\cite{Iwaniec2}, \cite{Kozachenko TVIMS}, \cite{Kozachenko4}, \cite{Ostrovsky comp sets}, \cite{Ostrovsky0}, \cite{Ostrovsky1},
\cite{Ostrovsky2}, \cite{Ostrovsky3}.\par

\vspace{4mm}

  \ Let us define the so - called {\it restricted Grand Lebesgue Space (rGLS) } $ \ GLS^{(S)}\psi. \ $
   Let $ \ S \ $ be some Borelian subset of whole semi-axis
$ \ [1,\infty): \ S \subset [1,\infty), \ $ containing
the initial point: $ \ \{1\} \in S. \ $ Define for arbitrary point $ \ p \ge 1 \ $

$$
p^+[S] = p^+[S](p) \stackrel{def}{=} \inf \{  t, \ t \in S, \ t \ge p \}
$$
and introduce the following variables

$$
Z[\psi,S] :=  \sup_{p \ge 1} \left[ \ \frac{\psi(p^+[S](p))}{\psi(p)} \ \right].
$$

\vspace{4mm}

 \ {\bf Definition 1.1.} \
The  restricted Grand Lebesgue Space (rGLS) $ \ G^{(S)}\psi \ $ consisted by definition on the (measurable) functions having a finite norm

$$
||f||G^{(S)}\psi \stackrel{def}{=} \sup_{p \in S} \left[ \ \frac{|f|_p}{\psi(p)} \ \right].
$$

\vspace{5mm}

 \ Evidently, these spaces are rearrangement invariant Banach functional spaces.\par

 \vspace{4mm}

\ {\bf Theorem 1.1.}  \  Suppose that  $ \ Z[\psi,S] < \infty. \ $   We propose that the classical GLS norm and introduced in
Definition 1.1 one are equivalent:

\vspace{3mm}

\begin{equation} \label{Zpsi}
 ||f||G^{(S)}\psi \le ||f||G\psi  \le Z[\psi,S] \ ||f||G^{(S)}\psi.
\end{equation}

\vspace{4mm}

 \ {\bf Proof.} \ The left - hand side of (\ref{Zpsi}) is obvious; we must ground opposite one. \par

  \ As long as the measure $ \ \mu \ $ is probabilistic, one can apply the Lyapunov's inequality

 $$
 |f|_p \le |f|_{p^+[S](p)}.
 $$
  \ Further, it follows from the definition of the value $ \ Z[\psi,S] \ $ for all the values $ \ p, \ p \ge 1 \ $

 $$
 \psi(p^+[S](p)) \le Z[\psi,S] \ \psi(p).
 $$
  \ We get substituting into the direct definition of the GLS norm

 $$
 \frac{|f|_p}{\psi(p)} \le \frac{|f|_{p^+[S]}}{\psi(p^+[S](p))/Z[\psi,S]} = Z[\psi,S] \cdot \frac{|f|_{p^+[S]}}{\psi(p^+[S](p))}.
 $$

  \ It remains to take supremum over  $ \ p \in [1,\infty). \ $ \par

 \vspace{5mm}

 \ We introduce now as a particular subspace of GLS  a so-called {\it discrete version} of Grand Lebesgue Space norm for the function $ \ \xi(\cdot). \ $  Let
 $ \ q = \{ q(1), q(2), q(3)\ldots  \}  \ $ be a strictly increasing sequence of real numbers such that

$$
q(1) = 1; \hspace{4mm} \lim_{m \to \infty} q(m) = \infty.
$$
  \ Denote by $ \ V = V( q) \ $ generating by $ \ q \ $ {\it partition } of semi - axis $ \ [1,\infty): \ $

$$
V = \{A(m)\}, \hspace{4mm} A(m) = [q(m), \ q(m+1)),
$$
so that

$$
m_1 \ne m_2 \ \Rightarrow    A(m_1) \cap A(m_2) = \emptyset
$$
and

$$
[1,\infty) = \cup_{m=1}^{\infty} A(m).
$$

\vspace{5mm}

\ {\bf Definition 1.2.} Define the following discrete Grand Lebesgue Space norm $ \ ||\xi||^q G\psi  \ $ for arbitrary r.v. $ \ \xi \ $  as a
particular case of rGLS

\vspace{3mm}

\begin{equation}  \label{Dis norm}
||\xi|| G^q\psi \stackrel{def}{=} \sup_{m=1,2,\ldots} \left\{ \ \frac{|\xi|_{q(m)}}{\psi(q(m))}  \ \right\}
\end{equation}
and correspondingly rearrangement invariant Banach functional {\it discrete Grand Lebesgue Space}

\begin{equation} \label{Dics space}
 G^q\psi \stackrel{def}{=} \{ \ \xi:  ||\xi|| G^q\psi < \infty. \ \}
\end{equation}

\vspace{5mm}

 \ For instance, if

$$
q=  q_0 = \{1,2,3,\ldots\},
$$
then

\begin{equation}  \label{Dis norm 0}
||\xi|| G^{q_0}\psi = \sup_{m=1,2,\ldots} \left\{ \ \frac{|\xi|_{m}}{\psi(m)}  \ \right\}.
\end{equation}

\vspace{4mm}

 \ Evidently, the introduced in (\ref{Dis norm}) functional is actually the rearrangement  invariant norm.\par

\vspace{4mm}

 \ Let us define the following function.

\begin{equation} \label{h funk}
h(x) = h[q,\psi](x) := \sup_m [ q(m) \ \ln x  - q(m) \ \ln \psi(q(m))], \ x \ge e.
\end{equation}
 \ It follows immediately from Tchebychev-Markov inequality that if $ \ 0 \ne \xi \in G^q\psi, \ $ then

\begin{equation} \label{tail estim}
{\bf P}(|\xi| \ge x) \le \exp \left( \ - h[q,\psi](x/||\xi||G^q\psi) \  \right), \ x \ge e \cdot ||\xi||G^q\psi.
\end{equation}
 \ The inverse proposition also holds true.  Indeed, if the r.v. $ \ \xi \ $ is such that

\begin{equation} \label{inverse tail estim}
\exists K = \const > 0 \ \Rightarrow {\bf P}(|\xi| \ge x) \le \exp \left( \ - h[q,\psi](x/K )\  \right), \ x \ge K \ e,
\end{equation}
then $ \ \xi \in   G\psi \ $ and moreover

$$
||\xi||G\psi \le C_1[\psi] \ \cdot K,
$$
see \cite{Kozachenko TVIMS}, \cite{Kozachenko4}.

\vspace{5mm}

\section{Main result.}

\vspace{5mm}

 \ We intend to establish in this section that under some natural conditions the ordinary and discrete Grand Lebesgue Spaces
 coincides up to norm equivalence. \par

  \ Note first of all that the inequality

\begin{equation} \label{triv est}
||f||G^q\psi \le ||f||G\psi
\end{equation}
is obvious. Let us derive an opposite inequality.\par
 \ Introduce a following variable.

$$
W = W[q,\psi]  \stackrel{def}{=} \sup_{m = 1,2,\ldots}  \left[ \  \frac{\psi(q(m+1))}{\psi(q(m))}   \ \right].
$$

 \ Obviously, $ \ W[q,\psi] > 1. \ $ \par

\vspace{4mm}

\ {\bf Theorem 2.1.} Suppose in addition that $ \ W[q,\psi] < \infty.   \ $ Then the ordinary  and discrete Grand Lebesgue Spaces norms are equivalent:

\begin{equation}  \label{equival}
||f||G^q\psi \le ||f||G\psi \le W[q,\psi] \cdot ||f||G^q\psi.
\end{equation}

\vspace{4mm}

 \ {\bf Proof.} It remains to ground the right-hand side of the proposition (\ref{equival}). Let $ \ W < \infty. \ $ We have

$$
||f||G\psi = \sup_{p \ge 1} \frac{|f|_p}{\psi(p)} = \sup_m \max_{p \in A(m)} \frac{|f|_p}{\psi(p)}.
$$

 \ One can apply the Lyapunov's inequality $ \  |f|_p \le |f|_{q(m+1)}, \ p \in A(m):  \ $

$$
||f||G\psi  \le \sup_m \max_{p \in A(m)} \frac{|f|_{q(m+1)}}{\psi(p)} =
$$

$$
\sup_m \left\{ \ |f|_{q(m+1)} \cdot \frac{1}{\min_{q \in A(m)} \psi(q) } \ \right\}
= \sup_m \left\{ \ |f|_{q(m+1)} \frac{1}{ \psi(q(m)) } \ \right\} \le
$$

$$
\sup_m \left\{ \ |f|_{q(m+1)} \cdot \frac{W[q,\psi]}{ \psi(q(m+1)) } \ \right\} \le
$$

$$
W[q,\psi] \ \sup_m \frac{|f|_{p(m+1)}}{\psi(p(m+1))} \le W[q,\psi] \ ||f||G^q\psi,
$$
Q.E.D.\par

\vspace{4mm}

 \ {\bf Remark 2.1.} \ The proposition of Theorem 2.1 remains true without condition of thee monotonicity of the
 generating function $\ \psi(\cdot). \ $ Namely, it is sufficient to suppose for certain sequence $ \ q = \{q(m)\} \ $ the condition
(\ref{positivity}) and to take instead the value $ \ W \ $ its modification

$$
\hat{W} = \hat{W}[q,\psi]  \stackrel{def}{=} \sup_{m = 1,2,\ldots}  \left[ \  \frac{\psi(q(m+1))}{\min_{p \in A(m)} \psi(p)} \ \right].
$$
so that

\begin{equation}  \label{new equival}
||f||G^q\psi \le ||f||G\psi \le \hat{W}[q,\psi] \cdot ||f||G^q\psi,
\end{equation}
if of course both the values $ \ \hat{W}[q,\psi] \ $ common with $ \ ||f||G^q\psi \ $  are finite.\par

\vspace{5mm}

\ {\bf Example 2.1.} The condition  $ \ \hat{W}[q,\psi]  < \infty \ $  is  satisfied if for instance

$$
\sup_m \frac{q(m+1)}{q(m)} < \infty
$$

and  the generating function $ \ \psi = \psi(p) \ $ has a popular form

$$
\psi(p) = p^{1/r} \ L(p), \ p \ge 1,
$$
where $ \ r = \const > 0, \ L = L(p) \ $ is some strictly positive: $ \ \forall p \ge 1 \ \Rightarrow \ L(p) > 0,  \ $ continuous slowly varying at infinity function.
For instance,

$$
L(p) = \ln^{\Delta}(2 + p), \ \Delta = \const.
$$

\vspace{3mm}

 \ One can take in this example

$$
q(m) = D^m -D + 1,
$$
for some fixed positive integer value $ \ D \ $ greatest or equal to 2:  $ \ D = 2,3,4,\ldots. \ $ \par

 \ The case $ \ L(p) = 1 \ $ and $ \ r = 2 \ $   correspondent to the  classical subgaussian random variable $ \ \xi. \ $ \par

\vspace{5mm}

\section{ The Restricted  Grand Lebesgue Spaces are also Banach algebras relative
to the convolution on unimodular locally compact groups.}

\vspace{5mm}

\ Let $ \ G \ $ be an unimodular local compact topological group equipped with bi-invariant (Borelian) Haar
measure $ \ \mu. \ $ Note that we do not suppose here that $ \ \mu \ $ is bounded. \par

\vspace{3mm}

 \ Define as ordinary the convolution between two measurable integrable functions $ \  f,g:G \to R \ $ as follows

$$
f*g(x) \stackrel{def}{=} \int_G f(y) g(y^{-1} \ x) \  d\mu(y).
$$

 \ The famous  Young (or Hausdorff-Young) inequality has the form

\begin{equation} \label{HausYoung}
|f*g|_r \le |f|_p \ |g|_q, \ 1 \le p,q,r \le \infty,
\end{equation}
where in addition

\begin{equation} \label{restiction}
 1+1/r = 1/p + 1/q,
\end{equation}
see e.g. \cite{Beckner}, \cite{Brascamp}, \cite{Fournier}, \cite{Hardy}, \cite{Zelazko}.\par

 \ It is proved in \cite{Formica Ban Alg}, that if the generating function of GLS $ \ \psi = \psi(p), \ 1 \le p <b, \ b = \const \in (1,\infty] $
 is continuous and such that $ \ \psi(1) = 1, \ $ then the Grand Lebesgue Space $ \ G\psi \ $ forms a {\it Banach algebra} relative the convolution, see in detail also
\cite{Gurkanli1}, \cite{Gurkanli2}, \cite{Gurkanli3}, \cite{Hewitt}. \par

 \vspace{3mm}

 \ {\it We intend in this section to extend this assertion into the restricted Grand Lebesgue Spaces } $ \ G^{(S)}\psi, \ $  defined before. \par

\vspace{4mm}

\ {\bf  Theorem 3.1. } Let the set  $ \ S \subset [1,\infty) \ $ be as before: $ \ \{1\} \in S \ $ and assume also $ \ \psi(1) = 1.\ $
Then the restricted Grand Lebesgue Space $ \ G^{(S)} \psi \ $ with correspondent norm

$$
||f||G^{(S)} \psi \stackrel{def}{=} \sup_{p \in S} \ \left\{ \ \frac{|f|_p}{\psi(p)} \ \right\}
$$
is also the Banach algebra relative the convolution operation on the
unimodular topological group  $ \ G \ $ equipped with the classical Haar's measure  $ \ \mu. \ $\par

\vspace{4mm}

\ {\bf Proof} is completely alike ones in \cite{Formica Ban Alg}. We must only ground the inequality

\begin{equation} \label{conv ineq}
||f*g|| G^{(S)} \psi \le ||f|| G^{(S)}\psi  \cdot ||g|| G^{(S)}\psi,
\end{equation}
see \cite{Zelazko}, chapter 2. \par
\ We can and will suppose without loss of generality

$$
||f|| G^{(S)} \psi = 1 =  ||g|| G^{(S)} \psi.
$$

 \ One can apply the inequality of Hausdorff-Young (\ref{HausYoung}):

$$
|f*g|_p \le |f|_p \cdot |g|_1, \ p \in S.
$$
 \ One has $ \ |g|_1 \le \psi(1) = 1, \ |f|_p \le \psi(p), \ $ therefore

$$
|f*g|_p \le \psi(p),
$$
following

$$
||f*g||G^{(S)} \psi \le 1 = ||f|| G^{(S)}\psi  \cdot ||g|| G^{(S)}\psi,
$$
Q.E.D. \par

\vspace{4mm}

\ {\bf Remark 3.1. } If the generating function  $ \ \psi = \psi(p) \ $ is not normalized: $ \ \psi(1) \ne 1, \ $ but is
such that

$$
\exists \lim_{ p \to 1+0} \psi(p) = \psi(1) \in (0,\infty),
$$
then the relation (\ref{conv ineq}) takes the form

\begin{equation} \label{not normalized}
||f*g||G^{(S)}\psi \le \psi(1) \cdot ||f||G^{(S)}\psi \cdot ||g||G^{(S)}\psi.
\end{equation}

\vspace{3mm}

 \ Herewith the condition $ \  \exists \ \psi(1+0) \in (0,\infty) \ $ is also necessary for the estimate of the form (\ref{not normalized})
 and the constant appears  therein $ \ \psi(1) \ $ is the best possible,  still in the particular case $ \ S = [1,\infty), \ $
 see \cite{Formica Ban Alg}. \par

\vspace{5mm}

\section{Concluding remarks.}

\vspace{5mm}

\hspace{6mm} {\bf A.} It is interest by our opinion to generalize obtained in this report  results
into the measurable sets having {\it an infinite measure,} as well as to consider the case when the support
of the generating function $ \ \psi(\cdot) \ $ is bounded:

$$
\psi(p) < \infty \ \Leftrightarrow 1 \le p < b, \ b = \const \in (1,\infty).
$$

\vspace{4mm}

 \ {\bf B.} It is interest  by our opinion to generalize also offered here notions and estimates into the
mixed (anisotropic) Lebesgue-Riesz spaces, to obtain for instance the CLT in these spaces, in the spirit
of an article \cite{Ostrovsky5}.\par

\vspace{6mm}

\vspace{0.5cm} \emph{Acknowledgement.} {\footnotesize The first
author has been partially supported by the Gruppo Nazionale per
l'Analisi Matematica, la Probabilit\`a e le loro Applicazioni
(GNAMPA) of the Istituto Nazionale di Alta Matematica (INdAM) and by
Universit\`a degli Studi di Napoli Parthenope through the project
\lq\lq sostegno alla Ricerca individuale\rq\rq (triennio 2015 - 2017)}.\par

\end{document}